\documentclass[12pt,twoside]{article}
\usepackage{times}
\usepackage{amsmath,amssymb}
\usepackage{amsthm}
\usepackage{mathrsfs}
\usepackage{color}
\usepackage{exscale}
\usepackage{relsize}
\usepackage{color}
\usepackage{fancyhdr}
\allowdisplaybreaks

\newcommand{\pf}{\noindent {\bf Proof. \hspace{2mm}}}
\newcommand{\ef}{ \hfill $ \Box $ \vskip 3mm}
\newcommand{\be}{\begin{equation}}
\newcommand{\ee}{\end{equation}}
\newcommand{\bea}{\begin{eqnarray}}
\newcommand{\eea}{\end{eqnarray}}

\newcommand{\bR}{{\mathbb R}}

\newcommand{\bS}{{\bf S}}

\def\nn{\nonumber}

\def\ve{\varepsilon}
\def\la{\lambda}
\def\La{\Lambda}

\def\q{\quad}
\def\qq{\qquad}

\def\Dl{\Delta}
\def\ve{\varepsilon}
\def\f{\frac}

\def\lt{\left}
\def\rt{\right}

\def\i{\infty}

\def\supp{\text{supp }}

\def\p{\partial}
\def\f{\frac}

\def\al{\alpha}

\def\o{\omega}

\def\lan{\langle}
\def\ran{\rangle}
\def\s{\sqrt}

 \allowdisplaybreaks

\begin{document}
 \footskip=0pt
 \footnotesep=2pt
\let\oldsection\section
\renewcommand\section{\setcounter{equation}{0}\oldsection}
\renewcommand\thesection{\arabic{section}}
\renewcommand\theequation{\thesection.\arabic{equation}}
\newtheorem{theorem}{\noindent Theorem}[section]
\newtheorem{lemma}{\noindent Lemma}[section]
\newtheorem*{claim}{\noindent Claim}
\newtheorem{proposition}{\noindent Proposition}[section]
\newtheorem*{inequality}{\noindent Weighted Inequality}
\newtheorem{definition}{\noindent Definition}[section]
\newtheorem{remark}{\noindent Remark}[section]
\newtheorem{corollary}{\noindent Corollary}[section]
\newtheorem{example}{\noindent Example}[section]
\newtheorem{conjecture}{\noindent Conjecture}[section]
\title{Blow up of solutions for semilinear wave equations with noneffective damping}

\author{Zijin Li $^{a,b,}$\footnote{E-mail:zijinli@smail.nju.edu.cn} ,\quad Xinghong Pan$^{c,}$\footnote{E-mail:xinghong{\_}87@nuaa.edu.cn}\vspace{0.5cm}\\
 \footnotesize $^a$Department of Mathematics and IMS, Nanjing University, Nanjing 210093, China.\\
 \footnotesize $^b$Department of Mathematics, University of California, Riverside, CA, 92521, USA.\\
\footnotesize $^c$Department of Mathematics, Nanjing University of Aeronautics and Astronautics, Nanjing 211106, China.
\vspace{0.5cm}
}

\date{}

\maketitle

\centerline {\bf Abstract} \vskip 0.3 true cm

In this paper, we study the finite-time blow up of solutions to the following semilinear
wave equation with time-dependent damping
\[
\p_t^2u-\Delta u+\frac{\mu}{1+t}\p_tu=|u|^p
\]
in $\mathbb{R}_{+}\times\mathbb{R}^n$. More precisely, for $0\leq\mu\leq 2,\mu \neq1$ and $n\geq 2$, there is no global solution for $1<p<p_S(n+\mu)$, where $p_S(k)$ is the $k$-dimensional Strauss exponent and a life-span of the blow up solution will be obtained. Our work is an extension of \cite{IS}, where the authors proved a similar blow up result with a larger range of $\mu$. However, we obtain a better life-span estimate when $\mu\in(0,1)\cup(1,2)$ by using a different method.


\vskip 0.3 true cm

\vskip 0.3 true cm


\section{Introduction}
\q In this paper, we deal with the following semilinear wave equation with time-dependent damping in multi-dimensions
\begin{equation}\label{Main}
\left\{
\begin{array}{ll}
\p_t^2u-\Delta u+\frac{\mu}{1+t}\p_tu=|u|^p,& (t,x)\in [0,\infty)\times\mathbb{R}^n,\\
\\
u(0,x)=\ve u_0(x),\,\,\p_tu(0,x)=\ve u_1(x),& x\in\mathbb{R}^n,\\
\end{array}
\right.
\end{equation}
where $n\geq2$, $\mu>0,\ \ve>0$. The initial data $(u_0(x),u_1(x))\in H^1(\bR^n)\times L^2(\bR^n)$ and supported in $\{x\in\bR^n||x|\leq R\}$.

Through this paper, we denote a generic constant by $C$ which may be different from line to line. $H^m(\mathbb{R}^n)$ denotes the usual Sobolev space with its norm
\be
\|f\|_{H^m}:=\sum_{|\al|=0}^m\|\p^\al_xf\|_{L^2}.\nn
\ee
We use $\|\cdot\|$ to denote $\|\cdot\|_{L^2}$ for convenience.

The critical exponent problem of \eqref{Main} for $p$ has been studied by many authors. Here the ``\emph{critical}'' means
there exists a critical exponent $p=p_c$ which satisfies:

{\it If $p>p_c$, the small-data solution of \eqref{Main} will exist globally in time; if $1<p\leq p_c$, the solution will blow up in finite time for data with positive average regardless of the smallness of the data.}

We note that for the linear part of \eqref{Main}, the damping term $\f{\mu}{1+t}$  is the borderline between the effective and noneffective dissipation, here effective means that the solution behaves like that of the corresponding parabolic equation and noneffective means that the solution behaves like that of the free wave equation. The asymptotic behavior of the solution relies heavily on the size of $\mu$.

Concretely, for the linear damped wave equation
\be
\p_{tt}u-\Dl u+\f{\mu}{(1+t)^{\beta}}u_t=0, \label{LN}
\ee
when $-1<\beta<1$, the asymptotic profile of the solution is given by the solution of the corresponding parabolic equation   $-\Dl v+\f{\mu}{(1+t)^{\beta}}v_t=0$. See \cite{Ma01}, \cite{Nk01}, \cite{Wj02} and references therein. On the other hand, when $\beta>1$, the asymptotic profile of the solution is the solution of the corresponding free wave equation $\p_{tt}v-\Delta v=0$. See \cite{Wj01}.

When $\beta=1$, the author in \cite{Wj03} proved that the solution of \eqref{LN} with data $(u,\p_tu)(0,x)=(u_0(x),u_1(x))$ satisfy
\bea
&&\|u(t)\|\leq C\Big(\|u_0\|+\|u_1\|_{H^{-1}}\Big)\lt\{
\begin{aligned}
&(1+t)^{1-\mu}\qq\ \mu\in(0,1),\\
&\ln(e+t)\qq\q \mu=1,\\
&1\qq \qq \qq \mu>1,
\end{aligned}
\rt.\nn\\
&& \|(\p_t,\nabla)u(t)\|\leq C\Big(\|u_0\|_{H^1}+\|u_1\|\Big)(1+t)^{\max\{-\f{\mu}{2},-1\}}.\nn
\eea

For the linear equation with space-dependent damping $\p_{tt}u-\Dl u+\f{\mu}{(1+|x|)^{1/2}}u_t=0$, readers can see \cite{ITY01} and reference therein for the asymptotic behavior of the solution.

Next we recall some results on the critical exponent problem to the semilinear wave equation $\p_{tt}u-\Dl u+\f{\mu}{(1+t)^\beta}\p_tu=|u|^p$, where $\mu>0$ .

 When $\beta=0$, Todorova-Yordanov \cite{TY} proved the critical exponent is given by $p_F=p_F(n):=1+2/n$. \cite{Zq} showed that the critical case $p=p_F$ also belonged to blow up case. \cite{IO} gave upper and lower life-span estimates of the solutions. Lin-Nishihara-Zhai in \cite{LNZ}(see also \cite{Nk02}) extend their results to the case $-1<\beta<1$ and proved that $p_F$ is still critical. When $\beta=1$, recently Wakasugi \cite{Wy01} proved that the critical exponnet is still $p_F$ for sufficiently large $\mu$. D'Abbicco \cite{Dm01} proved the global existence of small-data solutions under the assumption $\mu\geq n+2$ and $p>p_F$.  D'Abbicco-Lucente-Reissig \cite{DLR01} studied equations with more general effective damping
  \be
  \p_{tt}u-\Delta u+b(t)\p_tu=f(u).
  \ee
 See also \cite{FIW} \cite{II}. Ikeda and Wakasugi in \cite{IW} also considered the semilinear wave equation with time-dependent over damping.

  There are also many papers concerning about the critical exponent problem of the semilinnear wave equation with space-dependent or space-time dependent damping. See \cite{Wy02}, \cite{Wy03}, \cite{ITY02}, \cite{Lx}, \cite{IW01}, etc.

Here is our main theorem.

\begin{theorem}\label{blowup2}
 Let $\varphi(x)=\int_{\mathbb{S}^{n-1}}e^{x\cdot\omega}d\omega$. Suppose $0\leq\mu\leq2$, $\mu\neq1$, $n\geq2$, $(u_0(x), u_1(x))\in (H^1(\mathbb{R}^n), L^2(\mathbb{R}^n))$, both $u_0$ and $u_1$ supported in $\{x: |x|\leq R\}$, and
\be
\int_{\mathbb{R}^n}u_i(x)dx>0,\quad\int_{\mathbb{R}^n}u_i(x)\varphi(x)dx>0,\quad i=0,1.\nn
\ee
Then the solution of \eqref{Main} will blow up in finite time for any $\varepsilon>0$ when $1<p<p_S(n+\mu)$, where $p_S(\omega)$ is the positive root of the following equation with respect to $p$
\be
(\omega-1)p^2-(\omega+1)p-2=0.\nn
\ee
Moreover, the lifespan satisfies the following upperbound estimate
\be
T_\varepsilon\leq
C\varepsilon^{-\frac{2 p(p-1)  |1-\mu|}{2+( n + \mu +1 )p-(n + \mu-1 )p^2}},\quad 0\leq\mu\leq2,\q \mu\ne 1,\nn
\ee
where $C$ is a positive constant independent of $\varepsilon$.
\end{theorem}

\begin{remark}
Due to the fact $p_S(n+2\mu)<p_S(n+\mu)$ when $\mu>0$, our result partly improves the result by Lai-Takamura-Wakasa in \cite{LTW}, while their range of blow-up result for nonlinear exponent is $0<p<p_S(n+2\mu)$. Besides, when $n\leq 4$, our result is valid for a larger set of coefficient $\mu$ (except for $\mu=1$), since in \cite{LTW}, their blow up result is for $0<\mu<\frac{n^2+n+2}{2(n+2)}$.
\end{remark}

\begin{remark}
Not far from now, Ikeda-Sobajima in \cite{IS} proved a similar blow up result which contained the case $p=p_S(n+\mu)$, and they also obtained a larger range of $\mu$. However, our method is different from theirs and we achieve a better life-span estimate in the cases we consider.
\end{remark}

\begin{remark}
Note that $p_S(n+\mu)$ is a space shift of the Strauss exponent $p_S(n)$ and $p_S(n+\mu)>p_F(n)$ when $0<\mu<2,n\geq 2$. So \textbf{Theorem \ref{blowup2}} implies that our damping is no longer effective.  However, we still do not know whether $p_S(n+\mu)$ is the critical exponent of \eqref{Main}.
\end{remark}

In \cite{DLR02}, when $\mu=2$ and $n=2,3$, the authors have proved that the critical exponent is $p_c=p_S(n+2)$ which is a 2-dimension space shift of Strauss exponent $p_S$. Also D'Abbicco-Lucente \cite{DL} extend their results to any odd dimensions $n\geq 5$. However, when $0<\mu<2$, we can not say that $p_S(n+\mu)$ is the critical exponent of \eqref{Main}. But we can make a conjecture as follow which we think is reasonable in some sense.
\begin{conjecture}
When $0<\mu<2$, the critical exponent of \eqref{Main} is $p_c=p_S(n+\mu)$ which is the positive root of the following quadratic equation
\be
(n-1+\mu)p^2-(n+1+\mu)p-2=0.
\ee
\end{conjecture}

However, as far as the authors know, there is not effective way to prove the global existence of the solution when $p>p_S(n+\mu)$ by now.

The following is an outline of the proof of \textbf{Theorem \ref{blowup2}}. When $\mu\in[0,1)$, we make a time variable change $\La(t):=\f{ \langle t \rangle^{\ell+1}}{\ell+1}$ (we denote $\langle t\rangle:=1+t$ here and below), where $\ell=\f{\mu}{1-\mu}$ and set $w(t,x)=u(\La(t)-1,x)$. Then the Cauchy problem \eqref{Main} becomes a Cauchy problem of a semilinear free wave equation whose propagation speed is a polynomial of time
\be\label{change1}
\lt\{
\begin{aligned}
&\p^2_tw- \langle t\rangle^{2\ell}\Dl w= \langle t\rangle^{2\ell}|w|^p   \qq\qq\qq\qq\q\ (t,x)\in[t_0,\infty)\times\bR^{n}, \\
&w(t_0,x)= \varepsilon u_0(x), \ \p_t w(t_0,x)=\varepsilon(1-\mu)^{-\mu} u_1(x) \q x\in\bR^{n},
\end{aligned}
\rt.
\ee
where $t_0=(1-\mu)^{\mu-1}-1.$

When $\mu\in(1,2]$, similarly, we make the change of variables $\La(t):=\f{ \langle t \rangle^{\ell+1}}{\ell+1}$ and denote  $w(t,x):=\langle t\rangle u\big(\La(t)-1,x\big)$ with $\ell= \f{2-\mu}{\mu-1}$, once again the Cauchy problem \eqref{Main} becomes a Cauchy problem of a semilinear free wave equation whose propagation speed is a polynomial of time
\be\label{change2}
\lt\{
\begin{aligned}
&\p^2_tw- \langle t\rangle^{2\ell}\Dl w= \langle t\rangle^{2\ell-(p-1)}|w|^p   \qq\qq\qq\q\ (t,x)\in[t_0,\infty)\times\bR^{n}, \\
&w(t_0,x)= \varepsilon(\mu-1)^{1-\mu}\cdot u_0(x), \\
& \p_t w(t_0,x)=\varepsilon u_0(x)+\frac{\varepsilon u_1(x)}{\mu-1} \q x\in\bR^{n},
\end{aligned}
\rt.
\ee
where  $t_0=(\mu-1)^{1-\mu}-1$.

Using a key Lemma in He-Witt-Yin \cite{HWY} and the blow up trick originated from \cite{YZ} , we can prove that when $1<p<p_S(n+\mu)$, the solution of \eqref{change1}, \eqref{change2} will blow up in finite time which indicates that the solution of \eqref{Main} will also blow up for $1<p<p_S(n+\mu)$.

Our main theorem will be proved in the following section.

\section{Proof of the main theorem}

First we need the following ODE result.
\begin{lemma}\label{compare2}
(See \cite{St} and \cite{ZH}) Let $p>1$, $a\geq1$, and $(p-1)a>q-2$. If $f\in C^2([0,T))$ satisfies
\be
f(t)\geq\delta(t+R)^a, \label{4.1}
\ee
\be
f''(t)\geq m(t+R)^{-q}[f(t)]^p, \label{4.2}
\ee
where $\delta$, $R$ and $m$ are positive constants. Then $f(t)$ will blow up in finite time. The lifespan $T_\delta$ of $f(t)$ satisfies
\be
T_\delta\lesssim\delta^{-\frac{p-1}{(p-1)a-q+2}}.
\ee
\end{lemma}

\noindent \textbf{Proof of Theorem\ref{blowup2}:}

\q\  From the variable change in \eqref{change1} and \eqref{change2}, we only prove the case $0\leq\mu<1$ while the case $1<\mu\leq 2$ will be essentially the same.

In view of $\supp u_i\subseteq\ \{x\in\mathbb{R}^n\big||x|\leq R\}(i=0,1)$ and the finite propagation property of wave equations, we have for any fixed $t>0,$ the solution of \eqref{change1} is supported in $B(0,R+\La(t))$ where $\La(t)=\f{\lan t\ran^{\ell+1}}{\ell+1}$. We define
\be
F(t)=\int_{\bR^n}w(t,x)dx.
\ee
 Using integration by parts from $\eqref{change1}$ and H\"{o}lder inequality, we have
 \bea\label{est2.5}
F''(t)&=&\lan t\ran^{2\ell}\int_{\bR^n}|w|^pdx  \nn  \\
      &\geq&\lan t\ran^{2\ell}\f{|\int_{\bR^n}w(t,x)dx|^p}{\Big(\int_{|x|\leq R+\La(t)}dx\Big)^{p-1}} \nn \\
      &\geq&C\lan t\ran^{2\ell-n(\ell+1)(p-1)}|F(t)|^p,
 \eea
where we have used
\bea
|\int_{\bR^n}w(t,x)dx|&=&|\int_{\bR^n\cap \{|x|\leq R+\La(t)\}}w(t,x)dx| \nn \\
                      &\leq& C\Big(\int_{\bR^n}|w|^pdx\Big)^{\f{1}{p}}\Big(\int_{|x|\leq R+\La(t)}dx\Big)^{\f{p-1}{p}}.\nn
\eea
Next we introduce two test functions. The first one is the space test function
\be
\varphi(x)=\int_{\bS^{n-1}}e^{x\cdot\omega }d\o,
\ee
and it is well known that (see \cite{YZ}, section 2)
\be\label{ptest}
\varphi(x)\sim C_n|x|^{-\frac{(n-1)}{2}}e^{|x|},\quad\text{as }|x|\to+\infty;\quad\Delta\varphi=\varphi.
\ee
The second one is a modified Bessel function
\be
K_\al(t)=\int^\i_0e^{-t\cosh z}\cosh(\al z)dz,\q \al\in\bR,
\ee
which is a solution of the equation
\be
\Big(t^2\f{d^2}{dt^2}+t\f{d}{dt}-(t^2+\al^2)\Big)K_\al(t)=0.
\ee
From page 24 of \cite{EMOT}, when $\al>-\f{1}{2},$ we have
\be
K_\al(t)=\s{\f{\pi}{2t}}e^{-t}(1+O(t^{-1}))\qq \textmd{as} \q t\rightarrow\i.  \label{2.7}
\ee
We set
\be
\la(t)=C_\ell t^{\f{1}{2}}K_{\f{1}{2\ell+2}}\Big(\f{t^{\ell+1}}{\ell+1}\Big), \label{2.8}
\ee
where the constant $C_\ell$ is chosen such that $\la(t)$ satisfies
\be
\lt\{
\begin{aligned}
&&\la''(t)-t^{2\ell}\la(t)=0,\q t\geq t_0=(1-\mu)^{\mu-1}-1,  \\
&&\la(t_0)=1, \la(\i)=0.\q\q\q
\end{aligned}
\rt.\label{2.9}
\ee
We claim that $\la(t)$ has the following two properties (see Lemma 2.1 of \cite{HL}):\\
\indent 1. $\la(t)$, $-\la'(t)$ \emph{are both decreasing and} $\lim\limits_{t\rightarrow\i}\la(t)=\lim\limits_{t\rightarrow\i}\la'(t)=0,$ \\
\indent 2. \emph{There exists a constant} $C>1$ \emph{such that}
\be
\f{1}{C}\leq \f{-\la'(t)}{\la(t)t^{\ell}}\leq C,  \q \textmd{for} \q t\geq t_0.  \label{2.10}
\ee
Now we introduce the test function
\be\label{testfunction}
\psi(t,x)=\la(t)\varphi(x).
\ee
We have
\begin{lemma}\label{lem4.2}
Let $p>1$. $\psi$ is as in \eqref{testfunction}. Then there exists a $t_1>t_0$ such that for any $t\geq t_1$
\be\label{lem4.2e}
\int_{|x|\leq R+\La(t)}\psi(t,x)^{\f{p}{p-1}}dx\leq Ct^{-\f{\ell p}{2(p-1)}}(R+\La(t))^{n-1-\f{n-1}{2}\cdot\f{p}{p-1}}.
\ee
\end{lemma}
\pf By \eqref{2.7} and \eqref{2.8}, we have
\be
\la(t)\thicksim t^{-\f{\ell}{2}}e^{-\La(t)} \q \textmd{as} \ t\rightarrow \i.
\ee
This implies that there exists $t_1>0$ large enough, when $t>t_1$
\be\label{lem4.2e-1}
\begin{split}
\int_{|x|\leq R+\La(t)}\psi(t,x)^{\f{p}{p-1}}dx&=\la(t)^{\frac{p}{p-1}}\int_{|x|\leq R+\La(t)}\varphi(x)^{\f{p}{p-1}}dx\\
&\leq C t^{-\f{\ell p}{2(p-1)}}e^{-\frac{\La(t)p}{p-1}}\int_{|x|\leq R+\La(t)}\varphi(x)^{\f{p}{p-1}}dx\\
&:= C t^{-\f{\ell p}{2(p-1)}}e^{-\frac{\La(t)p}{p-1}} I_1.
\end{split}
\ee
As for $I_1$, by \eqref{ptest}, we have
\be\label{lem4.2e0}
\begin{split}
I_1=&\int_{|x|\leq R+\La(t)}\varphi(x)^{\f{p}{p-1}}dx\\
\leq & \,C\,\Big(\int^{\f{R+\La(t)}{2}}_0+\int^{R+\La(t)}_{\f{R+\La(t)}{2}}\Big)(1+r)^{n-1-\f{n-1}{2}\cdot\f{p}{p-1}}e^{\f{p}{p-1}r}dr :=I_{11}+I_{12}.
\end{split}
\ee
\be\label{lem4.2e1}
\begin{split}
I_{11}&\leq C(\La(t)+R)^{q_1}\int^{\f{R+\La(t)}{2}}_0e^{\f{p}{p-1}r}dr\\
&\leq C(\La(t)+R)^{q_1}\cdot e^{\frac{p}{2(p-1)}\La(t)},
\end{split}
\ee
where $q_1=\max\left(0, n-1-\f{n-1}{2}\cdot\f{p}{p-1}\right)$, and
\be\label{lem4.2e2}
\begin{split}
I_{12}&\leq C(\La(t)+R)^{n-1-\f{n-1}{2}\cdot\f{p}{p-1}}\int_{\f{R+\La(t)}{2}}^{R+\La(t)}e^{\f{p}{p-1}r}dr\\
&\leq C (\La(t)+R)^{n-1-\f{n-1}{2}\cdot\f{p}{p-1}}\cdot e^{\frac{p}{(p-1)}\La(t)}.
\end{split}
\ee
Combining \eqref{lem4.2e-1}, \eqref{lem4.2e0}, \eqref{lem4.2e1}, \eqref{lem4.2e2}, we finish our proof of of \textbf{Lemma \ref{lem4.2}}.

\ef
\noindent We denote
\be
F_1(t)=\int_{\bR^n}w(t,x)\psi(t,x)dx.  \label{2.11}
\ee
For the function $F_1(t),$ we have the following lower bound for $w(t,x)$ in both \eqref{change1} and \eqref{change2}.

\begin{lemma}
Under the assumption of \textbf{Theorem \ref{blowup2}}, there exists a $t_2>t_0$ such that
\be
F_1(t)\geq C\varepsilon t^{-\ell},\q t\geq t_2.   \label{2.12}
\ee
\end{lemma}

Readers can find the related proof in Lemma 2.3 of \cite{HWY}, and the $\varepsilon$ in \eqref{2.12} follows from a direct scaling since our initial data are of "$\varepsilon$-size". We omit the details here.
\ef
\noindent Using integration by parts, \eqref{2.12} and H\"{o}lder inequality, it follows that, for $t>t_2$
\bea\label{2.14}
F''(t)&=&\langle t\rangle^{2\ell}\int_{\bR^n}|w(t,x)|^pdx \nn  \\
      &\geq&t^{2\ell}\f{|F_1(t)|^p}{\big(\int_{|x|\leq R+\La(t)}\psi(t,x)^{\f{p}{p-1}}dx\big)^{p-1}} \nn  \\
      &\geq&C\varepsilon^p\f{t^{2\ell-p\ell}}{\big(\int_{|x|\leq R+\La(t)}\psi(t,x)^{\f{p}{p-1}}dx\big)^{p-1}},
\eea
where we have used
\bea
|F_1(t)|&=&\left|\int_{\bR^n\cap \{|x|\leq R+\La(t)\}}w(t,x)\psi(t,x) dx\right| \nn \\
        &\leq& \Big(\int_{\bR^n}|w|^p dx\Big)^{\f{1}{p}}\Big(\int_{|x|\leq R+\La(t)}\psi^{\f{p}{p-1}}dx\Big)^{\f{p-1}{p}}.\nn
\eea
Using \textbf{Lemma \ref{lem4.2}}, we get, for $t>t_3:=\max(t_1,t_2),$
\bea
F''(t)&\geq&C\varepsilon^pt^{2\ell-p\ell+\f{\ell}{2}p}(R+\La(t))^{-(n-1)(p-1)+\f{n-1}{2}p}  \nn  \\
     &=&C\varepsilon^pt^{2\ell-\f{\ell}{2}p}(R+\La(t))^{-(n-1)(p-1)+\f{n-1}{2}p}.    \label{2.17}
\eea
Integrating \eqref{2.17} twice from $t_3$ to $t$, and notice that $l=\frac{\mu}{\mu-1}$, we get
\be
F(t)\geq C\varepsilon^p(R+t)^{\frac{2+2n-p(n+\mu-1)}{2(1-\mu)}}+C_0(t-t_3)+C_1. \label{2.18}
\ee
When $1<p<p_S(n+\mu)$, $0\leq\mu<1$ and $n\geq 2$ it is easy to check that
\be
\frac{2+2n-p(n+\mu-1)}{2(1-\mu)}> 1,
\ee
then we have, for $t>t_3,$
\be
F(t)\geq  C\varepsilon^p(R+t)^{\frac{2+2n-p(n+\mu-1)}{2(1-\mu)}}.
\ee
So $F(t)$ satisfies \eqref{4.1} with
\be
a=\frac{2+2n-p(n+\mu-1)}{2(1-\mu)}.
\ee
In view of \eqref{est2.5}, we see $F(t)$ also satisfies \eqref{4.2} with
\be
q=-2\ell+n(\ell+1)(p-1)=\frac{2\mu+n-np}{\mu-1}.
\ee
Now the condition $a(p-1)>q-2$ means
\be
\frac{2+2n-p(n+\mu-1)}{2(1-\mu)}\cdot(p-1)>\frac{2\mu+n-np}{\mu-1}-2,
\ee
which is equivalent to
\be
(n+\mu-1)p^2-(n+\mu+1)p-2<0.
\ee

By using \textbf{Lemma \ref{compare2}}, we have the solution of problem \eqref{Main} with positive average data must blow up in finite time. Also the lifespan $T_\varepsilon$ satisfies
\be
\begin{split}
T_\varepsilon&\leq C(\varepsilon^p)^{-\frac{p-1}{(p-1)a-q+2}}\\
&\leq C \varepsilon^{-\frac{2 p(p-1)  (1 - \mu)}{2+( n + \mu +1 )p-(n + \mu-1 )p^2}}.\nn
\end{split}
\ee
When $1<\mu\leq 2$, we consider Cauchy problem \eqref{change2}. Similarly as \eqref{est2.5}, it follows that
\bea\label{est2.50}
F''(t)&=&\lan t\ran^{2\ell-(p-1)}\int_{\bR^n}|w|^pdx  \nn  \\
      &\geq&C\lan t\ran^{2\ell-(p-1)-n(\ell+1)(p-1)}|F(t)|^p\nn \\
      &=&C\lan t\ran^{\frac{3+n-\mu+(1-n-\mu)p}{\mu-1}}|F(t)|^p.
 \eea

Similarly as \eqref{2.18}, we have for $t>t_3$
\be
F(t)\geq C\varepsilon^p(R+t)^{\frac{2\mu+2n-(n+\mu-1)p}{2(\mu-1)}}+C_0(t-t_3)+C_1.
\ee
Again, when $1<p<p_S(n+\mu)$, $1<\mu\leq2$ and $n\geq 2$, it is easy to check
\be
\frac{2\mu+2n-(n+\mu-1)p}{2(\mu-1)}>1,
\ee
thus
\be
F(t)\geq C\varepsilon^p(R+t)^{\frac{2\mu+2n-(n+\mu-1)p}{2(\mu-1)}}, \quad\forall t>t_3.
\ee
Now the condition $a(p-1)>q-2$ in \textbf{Lemma \ref{compare2}} means
\be
\frac{2\mu+2n-(n+\mu-1)p}{2(\mu-1)}\cdot(p-1)>-\frac{3+n-\mu+(1-n-\mu)p}{\mu-1}-2,
\ee
which is equivalent to
\be
(n+\mu-1)^2p-(n+\mu+1)p-2<0.
\ee
Then \textbf{Lemma \ref{compare2}} leads to the blow up result and the following estimate of lifespan
\be
T_\varepsilon\leq C\varepsilon^{-\frac{2 p(p-1)  (\mu-1)}{2+( n + \mu +1 )p-(n + \mu-1 )p^2}}.\nn
\ee
We complete the proof of \textbf{Theorem\ref{blowup2}}.

\rightline{$\Box$}

\indent

{\bf Acknowledgement.}

The authors wish to thank Professor Huicheng Yin in Nanjing Normal University for his constant encouragement on this topic. They are also grateful to the department of mathematics UC Riverside for the hospitality during their visit, when part of the paper was written.

Zijin Li is supported by China Scholarship Council (File No. 201606190089). 

\
\\

\end{document}